\documentclass[11pt, a4paper]{article}
\usepackage[utf8]{inputenc}
\usepackage{bbm, bm}

\usepackage[margin = 1in]{geometry}
\usepackage{amsfonts,amsmath,amssymb,amsthm}
\usepackage{mathtools}
\usepackage{mathrsfs}
\usepackage{graphicx,subcaption}
\usepackage{booktabs}
\usepackage{tabularx}
\usepackage{lipsum}
\usepackage{hyperref}
\usepackage{multicol}
\hypersetup{
   colorlinks = true,
   citecolor = blue,
   urlcolor = purple,
   linkcolor = red
}

\usepackage{mathrsfs}
\usepackage{multirow}
\usepackage[makeroom]{cancel}
\usepackage{tikz}
\usepackage{framed,color}
\definecolor{shadecolor}{rgb}{1,0.8,0.3}
\usepackage{fancybox}
\usepackage{caption}
\usepackage{algorithm,algorithmic} 
\usepackage{multicol}
\usepackage{aliascnt}
\usepackage{setspace}
\usepackage{longtable}
\allowdisplaybreaks
\usepackage{comment}

\usetikzlibrary{decorations.markings,arrows}

\allowdisplaybreaks

\title{\textbf{Some Asymptotic Results on Multiple Testing under Weak Dependence}}

\date{}
\vspace{5mm}
\usepackage{authblk}
\author[1]{Swarnadeep Datta\footnote{\href{swarnadeepdatta0122@gmail.com}{swarnadeepdatta0122@gmail.com}}}
\author[2]{Monitirtha Dey\footnote{\href{monitirthadey3@gmail.com}{monitirthadey3@gmail.com}, \href{mdey@uni-bremen.de}{mdey@uni-bremen.de}}}
\affil[1]{\footnotesize Interdisciplinary Statistical Research Unit, Indian Statistical Institute, Kolkata, India}
\affil[2]{\small Institute for Statistics, University of Bremen, Bremen, Germany}
 
\usepackage{natbib}
\setcitestyle{authoryear, open={(},close={)}}

\begin{document}

\maketitle
\theoremstyle{plain}
\newtheorem{axiom}{Axiom}
\newtheorem{remark}{Remark}
\newtheorem{corollary}{Corollary}[section]
\newtheorem{claim}[axiom]{Claim}
\newtheorem{theorem}{Theorem}[section]
\newtheorem{lemma}{Lemma}[section]
\newtheorem{test}{Test Procedure}
\newtheorem{proposition}{Proposition}

\newaliascnt{lemmaa}{theorem}
\newtheorem{lemmaa}[lemmaa]{Theorem}
\aliascntresetthe{lemmaa}
\providecommand*{\lemmaautorefname}{Lemma}
\providecommand*{\corollaryautorefname}{Corollary}
\providecommand*{\testautorefname}{Test Procedure}
\providecommand*{\theoremautorefname}{Theorem}
\providecommand*{\propositionautorefname}{Proposition}
\providecommand*{\remarkautorefname}{Remark}



\theoremstyle{plain}
\newtheorem{exa}{Example}
\newtheorem{rem}{Remark}

\theoremstyle{definition}
\newtheorem{definition}{Definition}
\newtheorem{example}{Example}

\begin{abstract}
This paper studies the means-testing problem under weakly correlated Normal setups. Although quite common in genomic applications, test procedures having exact FWER control under such dependence structures are nonexistent. We explore the asymptotic behaviors of the classical Bonferroni (when adjusted suitably) and the Sidak procedure; and show that both of these control FWER at the desired level exactly as the number of hypotheses approaches infinity. We derive analogous limiting results on the generalized family-wise error rate and power. Simulation studies depict the asymptotic exactness of the procedures empirically.
\end{abstract}

\section{Introduction}

Large-scale multiple testing problems under dependence are a staple of modern statistical science \citep{Dickhaus2026, FDR2007, FDR2009}. Dependence among observations or test statistics has been a springboard for the theoretical and applied developments of simultaneous statistical inference \citep{DickhausThesis, Dickhaus}. Often, there exists a weak dependence among the observations of interest. \cite{Proschan} mention that the correlation between two single nucleotide polymorphisms (SNP) is generally believed to decrease with genomic distance. Several authors 
have argued that for the large numbers of markers typically 
used for a genome-wide association study (GWAS), the test statistics are weakly correlated due to this largely local presence of correlation between SNPs \citep{DeyBhandari2025gfwer, StoreyTib}. 

There have been a few works on weakly dependent structures with respect to asymptotic false discovery rate (FDR) control. \cite{Farcomeni} shows that a certain degree of dependence is allowed among the test statistics when the number of tests is large, without any need for a correction to the traditional procedures. \cite{Storey2004} show that
 linear step-up (LSU) and  plug-in LSU tests control
the FDR asymptotically under weak dependence, assuming that the proportion of rejected null hypotheses is asymptotically larger than zero in some sense. \cite{GONTSCHARUK} show that weak dependence is not sufficient for FDR
control if the proportion of rejected nulls converges to zero with positive probability. However, similar theoretical results on asymptotic familywise error rate (FWER) control are scarce. \cite{DasBhandari2025} explore how closely the FWER of Bonferroni method resembles its behavior under independence and introduces an asymptotic correction factor to improve accuracy in nearly independent cases. 

A series of recent works viz. \cite{deybhandari, deybhandaristpa, DeyBhandari2025gfwer, deystpa} has elucidated that the Bonferroni method and the class of stepwise multiple testing procedures  (MTP) have asymptotically zero family-wise error rate (FWER) under correlated Normal scenarios where the infimum of correlations is considered to be strictly positive. However, their limiting results do not help in devising an MTP that is asymptotically exact under equicorrelation, i.e., which has limiting FWER exactly equal to $\alpha$ under equicorrelation. Inspired by such problems, the authors propose a simple single-step MTP that asymptotically controls the FWER at the desired level exactly under the equicorrelated multivariate Gaussian setup in \cite{DATTA2026110609}. However, such results guaranteeing asymptotic exact FWER control are nonexistent for weakly dependent scenarios. Also, asymptotic behaviors of Bonferroni-type procedures have not yet been studied under general weak dependence structures. This work fills both these gaps by elucidating that both the Sidak procedure and the classical Bonferroni method (with a little adjustment) have asymptotic FWER exactly equal to $\alpha$ under weakly dependent scenarios.

This work is organized as follows. Section \ref{sec2} introduces the testing framework formally. Section \ref{sec3} presents the main results on asymptotically exact FWER control in the one-sided setting. Section \ref{sec4} investigates the power properties. Section \ref{sec5} studies the same problem in the two-sided setting. Section \ref{sec6} depicts the empirical performance of the procedures through simulations. Section \ref{sec7} concludes the paper. The proofs are deferred to the appendix. 

\section{Preliminaries \label{sec2}}
Throughout this paper, $\phi$ and $\Phi$ denote the p.d.f and c.d.f of $N(0,\,1)$ distribution respectively. Let $\mathcal{I}$ denote the set $\{1,\,2,\,\dots,\,n\}$. Consider $n$ observations
$$X_i \sim N(\mu_i, 1),   i \in \mathcal{I}.$$
The elements of the covariance matrix $\Sigma_n$ are given by 
\begin{align*}
& Corr(X_i,\,X_j)=\rho_{ij}, \text{ with } \rho_{ij}\in(-1,\,1) \quad \forall\  i\ne j.
\end{align*}
Let $\displaystyle\rho_m=\sup_{1\leq i\leq n-m}|\rho_{i\, \overline{i+m}}|$ and also $\displaystyle\gamma=\sup_{n\ge 1}\rho_n<1$. We assume that $\Sigma_n$  satisfies the following \textit{weak dependence condition}:
\begin{equation}\label{wd_cond}\rho_m=o \left(\frac{1}{\log m}\right) \quad \forall\ 1\leq m \leq n \hspace{2mm} \text{as $n$ grows.} \end{equation}
In sections \ref{sec3} and \ref{sec4}, we consider the one-sided testing problem
$$H_{0i} : \mu_i = 0 \quad \text{vs} \quad H_{1i} : \mu_i > 0, \quad i \in \mathcal{I}.$$
We focus on the corresponding both-sided problem 
$$H_{0i} : \mu_i = 0 \quad \text{vs} \quad H_{1i} : \mu_i \neq 0, \quad i \in \mathcal{I}$$
in section \ref{sec6}. In both settings, the global null $H_0 = \bigcap_{i=1}^{n} H_{0i}$ asserts that each $\mu_i$ is zero.

This paper aims to obtain valid MTPs for which FWER converges to the desired level, under the weak dependence structure depicted in \eqref{wd_cond}. 

\section{Main Results \label{sec3}} Throughout this work, let $\mathcal{I}_0$ and $\mathcal{I}_1$ respectively denote the set of indices of true and false nulls. Let $n_0$ and $n_1(=n-n_0)$ be the cardinalities of these two sets, respectively. Also, unless otherwise mentioned, by the phrase \textit{weakly dependent Gaussian sequence}, we mean any Gaussian sequence satisfying \eqref{wd_cond}. We consider the following two testing procedures:

\begin{test}[Adjusted Bonferroni]\label{test_bon} For  each $i\in\mathcal{I}$, reject $H_{0i}$ if $X_i>c_{Bon}(n,\,\alpha)\coloneqq\Phi^{-1}(1-\frac{-\log\,(1-\alpha)}{n})$. 
\end{test}

\begin{test}[Sidak]\label{test_sidak}
 For  each $i\in\mathcal{I}$, reject $H_{0i}$ if $X_i>c_{Sid}(n,\,\alpha)\coloneqq\Phi^{-1}((1-\alpha)^{1/n})$.
\end{test}

\noindent We first provide the following result on the asymptotic distribution of the $k$-th largest value for a weakly dependent standard Gaussian sequence: 
\begin{theorem}\label{probabilistic_k_dn}
    Let $\left\{X_n\right\}$ be a weakly dependent standard Gaussian sequence. Suppose $\left\{u_n\right\}$ be a sequence and $\left\{d_n\right\}$ be a positive integer sequence of n such that $d_n\to\infty$ and ${d_n}(1-\Phi(u_n))\to \tau\text{ as }n\to\infty$ for some $\tau \in [0,\,\infty]$. For $1 \leq k \leq d_n$, let $M_{d_n}^k$ be the $k$-th largest among $X_1,\,X_2,\,\dots,\,X_{d_n}$. Then, for each fixed $k$, we have
\begin{equation}
\displaystyle\lim_{n\to\infty}\mathbb{P}\left\{M_{d_n}^k \leq u_n\right\} =\begin{cases}
    e^{-\tau} \sum_{s=0}^{k-1} \frac{\tau^s}{s!}\quad&\text{if }\tau\in(0,\,\infty),\\
    1\quad&\text{if }\tau=0,\\
    0\quad&\text{if }\tau=\infty.\\
\end{cases}
\end{equation}
\end{theorem}

Let $i_1<i_2<\cdots<i_{n_0}$ be the indices in $\mathcal{I}_0$. Then, we re-notate $X_i$'s as $Z_j=X_{i_j}$, $j=1,\,2,\,\dots,\,n_0$. The weak dependence structure gives $Corr(Z_i,\,Z_j)\leq\rho_{|i-j|}$. Hence, as $n$ (and thereby $n_0$) grows, $Corr(Z_1,\,Z_{m})\leq\rho_{m-1}=o(\frac{1}{\log (m-1)})$ satisfying the weak dependence structure for $\{Z_m\}$ as well.

\noindent Now, FWER of any procedure with common right-sided cutoff $\tau_n$ can be expressed as:
\begin{align}
FWER(n, \tau_n, \alpha, \Sigma_n)&=\mathbb{P} \left(\displaystyle\bigcup_{i=1}^{n_0}\{Z_i>\tau_n\} \right)=1-\mathbb{P}\left(\displaystyle\bigcap_{i=1}^{n_0}\{Z_i\leq \tau_n\}\right)\notag\\
&=1-\mathbb{P}\left(M_{n_0}\leq \tau_n\right), \text{ where }M_{n_0}= \max_{1 \leq i \leq n_0} Z_i.
\end{align}

\noindent Note that for any $\alpha\in(0,\,1)$,
\begin{align}\label{sid_exp}
    (1-\alpha)^{\frac{1}{n}}=\exp(\frac{\log(1-\alpha)}{n})&=1-\frac{-\log(1-\alpha)}{n}+O(\frac{1}{n^2})\notag \\
    &=1-\frac{-\log(1-\alpha)-O(\frac{1}{n})}{n}.
\end{align}

 \noindent Hence, we can express both $c_{Bon}(n,\,\alpha)$ and $c_{Sid}(n,\,\alpha)$ as $c_n(\alpha)=\Phi^{-1}(1-\frac{-\log(1-\alpha)-o(1)}{n})$.

We note that $n\left(1-\Phi(c_n(\alpha))\right)\longrightarrow-\log(1-\alpha)\text{ as }n\to\infty$. Taking $u_n=c_n(\alpha)$, $d_n=n_0$ and $k=1$ in \autoref{probabilistic_k_dn}, we obtain the quintessential result of this work.

\begin{theorem}\label{fwer_bon_2}
    Under the weakly dependent standard Gaussian setting with correlation matrix $\Sigma_n$, both the adjusted Bonferroni procedure(\autoref{test_bon}) and the Sidak procedure(\autoref{test_sidak}) are asymptotically exact, i.e.,
    $$\displaystyle\lim_{n\to\infty}\,FWER(n, c_{Bon}, \alpha,\,\Sigma_n)=\lim_{n\to\infty}\,FWER(n, c_{Sid}, \alpha,\,\Sigma_n)=\alpha,$$
    under any configuration of true and false null hypotheses satisfying $\lim_{n\to\infty}n_0/n=1$.
\end{theorem}

\begin{remark}\label{adj_tests}
Let $p_0:= \lim_{n\to\infty}n_0/n$. Notably, the FWER of \autoref{test_bon} for this weak dependence setup, under any configuration of true and false null hypotheses, does not necessarily converge to target $\alpha$ for any $p_0>0$. \autoref{fwer_bon_2} states that this convergence is valid for $p_0=1$. If $p_0$ is known a priori, then one may suitably tweak the common cutoffs as mentioned in \autoref{test_bon} and \autoref{test_sidak} to guarantee exact convergence. Indeed, the procedures utilizing the cutoffs $c_{Bon}(\alpha,\,p_0)\coloneqq\Phi^{-1}(1-\frac{-\log\,(1-\alpha)}{np_0})$ and $c_{Sid}(\alpha,\,p_0)\coloneqq\Phi^{-1}((1-\alpha)^{1/np_0})$ have their FWERs converging to $\alpha$.
\end{remark}

\begin{remark}\label{rate_of_convergence}
    Let $\nu$ be some constant such that $0<\nu<\frac{1-\gamma}{1+\gamma}$, $\gamma=\sup_{n \geq 1}\rho_{n}<1$ and also let $\gamma_n = \sup_{m \geq n} \rho_{m}$. Then the rate of convergence in \autoref{fwer_bon_2} can be given as
    $$|FWER(n, c_n(\alpha),\alpha, \Sigma_n) -\alpha| \leq l\cdot R_n\quad \text{for some $l>0$ as $n\to \infty$}$$
    where
    $$R_n=\max\left\{n^{\frac{1+\nu-2}{1+\gamma}}(\log n)^{\frac{1}{1+\gamma}},\ \gamma_{\left[n^\nu\right]}\log n^{\nu},\ 1-\frac{n_0}{n},\ \frac{1}{n}\right\}.$$
 Here $c_n(\alpha)\in \{ c_{Bon}(n,\,\alpha),\,c_{Sid}(n,\,\alpha)\}$.
\end{remark}

\begin{remark}\label{k_fwer}
    For the Lehmann-Romano procedure \citep{Lehmann2005},
\begin{align*}k\text{-FWER}(n, \alpha, \Sigma_n) &=\mathbb{P}_{\Sigma_n}\left(X_{i}>\Phi ^{-1}(1-k\alpha/n) \hspace{2mm} \text{for at least $k$}
\hspace{2mm}\text{$i$'s}\in\mathcal{I}_0\right)\\
&= 1 - \mathbb{P}_{\Sigma_n}\left(M_{n_0}^k \leq \Phi^{-1}(1-k\alpha/n)\right)\\
&\quad \quad \quad [M_{n_0}^k:k\text{'th} \text{ maximum among $Z_1,\,\dots,\,Z_{n_0}$.}]\end{align*}
Considering $u_n=\Phi ^{-1}(1-k\alpha/n)$ and $d_n=n_0$ in \autoref{probabilistic_k_dn} directly gives the asymptotic $k$-FWER control for the Lehmann-Romano procedure under the weakly dependent standard Gaussian setting with correlation matrix $\Sigma_n$, i.e.,
  \[k\text{-FWER}(n, \alpha, \Sigma_n)  \longrightarrow \displaystyle1- e^{-k\alpha} \sum_{s=0}^{k-1} \frac{(k\alpha)^s}{s!} \text { as } n \rightarrow \infty .\]
\end{remark}

\section{Power Analysis \label{sec4}}

The simultaneous inference literature has several notions of power \citep{Dudoit}. In this work, we shall work with \textit{AnyPwr}, which is defined as the probability of making at least one true rejection. We have $X_i\sim N(\mu_i, 1)$, $\mu_i>0$ for $i\in\mathcal{I}_1$, where $\mathcal{I}_1$ denotes the set of indices of the originally false nulls. Hence, $|\mathcal{I}_1|=n_1$.
The following result describes the asymptotic powers of \autoref{test_bon} and \autoref{test_sidak}.

\begin{theorem}\label{pwr_thm_1}
    Suppose $n_1\to\infty$ and $\frac{n_1}{n}\to p_1\in(0,\,1]$ as $n\to\infty$. Then, under the weakly correlated Gaussian setting, for both \autoref{test_bon} and \autoref{test_sidak}, one has $\displaystyle\lim_{n\to\infty} AnyPwr=1$ if $\displaystyle\lim_{n_1\to\infty}\,\frac{\sqrt{2\,\log n_1}}{\mu_{n_1}}<1$ where $\mu_{n_1}=\max\{\mu_i,\ i\in\mathcal{I}_1\}$.
\end{theorem}




We discuss now the asymptotic of power for a more general setup of the non-null means. The following result would be crucial for that. 

\begin{proposition}\label{Phi^dn}
Suppose, for any $\beta>0$, $\beta_n=\beta+o(1)>0$ and $c_{\beta_n,\,n} = \Phi^{-1}(1-\beta_n/n)$. Also suppose $\{d_n\}_{n\geq1}$ is a sequence of $n$ which diverges to $\infty$ as $n \to \infty$. Then, if for some $t\in\mathbb{R}$, $\displaystyle\lim_{n\to\infty}\frac{d_n}{n}\cdot e^{-t\sqrt{2\log n}}=\infty$, one has
\[\displaystyle\lim_{n\to\infty}d_n \cdot\left(1-\Phi\left(c_{\beta,\,n}+t\right)\right)=\infty\]
\end{proposition} 

Let $i_1<i_2<\cdots<i_{n_1}$ be the indices in $\mathcal{I}_1$. Then we re-notate $X_i$'s as $Y_j=X_{i_j}$, and also, throughout this section, we would use $\mu_j$ to denote $\mu_{i_j}$, $j=1,\,2,\,\dots,\,n_1$. Evidently, $Corr(Y_i,\,Y_j)\leq\rho_{|i-j|}$ and as $n$ and thereby $n_1$ grows, $Corr(Y_1,\,Y_{m})\leq\rho_{m-1}=o(\frac{1}{\log (m-1)})$.

One may write $Y_j=Z_j+\mu_j,\ Z_j\sim N(0,\,1)\ \forall\ j=1,\,2,\,\dots,\,n_1$. Clearly, $Corr(Y_i,\,Y_j)=Corr(Z_i,\,Z_j)$ $\forall\ i, j=1,\,2,\,\dots,\,n_1$.

Now, for any testing procedure with common right-sided cutoff $\tau_n$, we have the following inequality.
\begin{align}\label{pwr_eq_2}
    1-AnyPwr&=\mathbb{P}(\bigcap_{i=1}^{n_1}\{Z_i\leq \tau_n-\mu_i\})\nonumber\\
    &\leq\mathbb{P}(\bigcap_{i=1}^{n_1}\{Z_i\leq \tau_n-\mu\})\quad\text{[assuming $\displaystyle \mu=\lim_{n\to\infty}\min_{i\in\mathcal{I}_1}\mu_i>0$]}\nonumber\\
    &=\mathbb{P}(Z_{(n_1)}\leq \tau_n-\mu)\quad\text{[where $Z_{(n_1)}=\displaystyle\max_{i=1,\,2,\,\dots,\,n_1}Z_i$]}.
\end{align}
We put $d_n=n_1$, $t=-\mu$ and $\beta_n=-\log(1-\alpha)-o(1)$. \autoref{Phi^dn} implies that
$$\displaystyle\lim_{n\to\infty}d_n \cdot\left(1-\Phi\left(c_{\beta,\,n}+t\right)\right)=\infty.$$
Applying \autoref{probabilistic_k_dn} on $u_n=c_n(\alpha)-\mu$ and $k=1$, we get the following result viding \eqref{pwr_eq_2}.

\begin{theorem}\label{pwr_thm_2}
    Suppose $\displaystyle\lim_{n\to\infty}\min_{i\in\mathcal{I}_1}\mu_i=\mu>0$. Then, for both \autoref{test_bon} and \autoref{test_sidak}, one has $\displaystyle\lim_{n\to\infty}\,AnyPwr=1$ if $\displaystyle\lim_{n\to\infty}\frac{n_1}{n}e^{\mu\sqrt{2\log n}}=\infty$.
    \end{theorem}

\section{Extension to Both-sided Testing \label{sec5}}

Up to this point, we have considered only the one-sided testing problem, i.e., $H_{1i} : \mu_i > 0,\ i \in \mathbb{I}$. However, one often encounters both-sided testing situations:
$$H_{0i} : \mu_i = 0 \quad \text{vs} \quad H_{1i} : \mu_i \neq 0, \quad i \in \mathcal{I}.$$
We shall denote the FWER, $k$-FWER, and AnyPwr in this setting by $FWER_{BS}$, $k\text{-}FWER_{BS}$ and $AnyPwr_{BS}$ respectively.
For this problem, we consider the corresponding testing procedures.

\begin{test}[Adjusted Bonferroni]\label{test_bon_bs}
    For any covariance matrix (or correlation matrix) $\Sigma_n$ satisfying the weak dependence condition \eqref{wd_cond}, taking the common cutoff $\{c_{Bon}(n,\,\alpha)\}_{n\geq 1}$ as defined in \autoref{test_bon}, we reject $H_{0i}$ if $| X_i|>c_{Bon}(2n, \alpha)$ for each $i\in\{1,\,2,\,\dots,\,n\}$.
\end{test}

\begin{test}[Sidak]\label{test_sidak_bs}
    For any covariance matrix (or correlation matrix) $\Sigma_n$ satisfying the weak dependence condition \eqref{wd_cond}, taking the common cutoff $\{c_{Sid}(n,\,\alpha)\}_{n\geq 1}$ as defined in \autoref{test_sidak}, we reject $H_{0i}$ if $| X_i|>c_{Sid}(2n, \alpha)$ for each $i\in\{1,\,2,\,\dots,\,n\}$.
\end{test}

Analogous to section \ref{sec3}, we first provide the result on asymptotic distribution of the $k$-th largest among the absolute values of a weakly dependent standard Gaussian sequence.

\begin{theorem}\label{probabilistic_k_dn_bs}
    Consider the assumptions as in \autoref{probabilistic_k_dn}. Suppose $L_{d_n}^k$ be the $k$-th largest among $|X_1|,\,|X_2|,\,\dots,\,|X_{d_n}|$ for $1 \leq k \leq d_n$. For each fixed $k$, we have
\begin{equation}
\displaystyle\lim_{n\to\infty}\mathbb{P}\left\{L_{d_n}^k \leq u_n\right\} =\begin{cases}
    e^{-2\tau} \sum_{s=0}^{k-1} \frac{\left(2\tau\right)^s}{s!}\quad&\text{if }\tau\in(0,\,\infty),\\
    1\quad&\text{if }\tau=0,\\
    0\quad&\text{if }\tau=\infty.\\
\end{cases}
\end{equation}
\end{theorem}
    
Consider the rearrangements and re-notations of $X_i$'s corresponding to the originally ture nulls mentioned in section \ref{sec3}. The FWER of any procedure with common left-sided and right-sided cutoffs $-\tau_n$ and $\tau_n$, respectively, can be expressed in the following form:
\begin{align}\label{eq_fwer_ex}
FWER_{BS}(n, \tau_n, \alpha, \Sigma_n)&=\mathbb{P} \left(\displaystyle\bigcup_{i=1}^{n_0}\{|Z_i|>\tau_n\} \right)=1-\mathbb{P}\left(\displaystyle\bigcap_{i=1}^{n_0}\{|Z_i|\leq \tau_n\}\right)\notag\\
&=1-\mathbb{P}\left(L_{n_0}\leq \tau_n\right), \text{ where }L_{n_0}= \max_{1 \leq i \leq n_0} |Z_i|.
\end{align}

Then, taking $u_n=c_{2n}(\alpha)=\Phi^{-1}(1-\frac{-\log(1-\alpha)-o(1)}{2n})$, $d_n=n_0$ and $k=1$ in \autoref{probabilistic_k_dn_bs}, we get the following result for the two-sided testing problem.

\begin{theorem}\label{bs_thm_1}
    For both \autoref{test_bon_bs} and \autoref{test_sidak_bs} for the both-sided problem, under general correlated Gaussian setting with $\Sigma_n$ satisfying the weak dependence condition \ref{wd_cond}, $$\displaystyle\lim_{n\to\infty}FWER_{BS}(n, c_{Bon}, \alpha, \Sigma_n)=\lim_{n\to\infty}FWER_{BS}(n, c_{Sid}, \alpha, \Sigma_n)=\alpha,$$
    under any configuration of true and false null hypotheses for which $\lim_{n\to\infty}n_0/n=1$.
\end{theorem}

\begin{remark}\label{rate_of_convergence_bs}
    Let $R_n$ be as defined in \autoref{rate_of_convergence}. Then the rate of convergence in \autoref{bs_thm_1} is also given as
        $$|FWER_{BS}(n, c_n(\alpha),\alpha, \Sigma_n) -\alpha| \leq l\cdot R_n\quad \text{for some $l>0$ as $n\to \infty$}.$$
\end{remark}

\begin{remark}\label{k_fwer_bs}
    For the Lehmann-Romano procedure \citep{Lehmann2005} for this two-sided testing problem,
\begin{align*}k\text{-FWER}_{BS}(n, \alpha, \Sigma_n) &=\mathbb{P}_{\Sigma_n}\left(|X_{i}|>\Phi ^{-1}(1-k\alpha/2n) \hspace{2mm} \text{for at least $k$}
\hspace{2mm}\text{$i$'s}\in\mathcal{I}_0\right)\\
&= 1 - \mathbb{P}_{\Sigma_n}\left(L_{n_0}^k \leq \Phi^{-1}(1-k\alpha/2n)\right)\\
&\qquad\Bigl[
L_{n_0}^k:\ k\text{'th maximum among } |Z_1|,\,\dots,\,|Z_{n_0}|\text{as defined in Section~\ref{sec3}.}
\Bigr]\end{align*}
Considering $u_n=\Phi ^{-1}(1-k\alpha/2n)$ and $d_n=n_0$ in \autoref{probabilistic_k_dn_bs} directly gives the asymptotic $k$-FWER control for the Lehmann-Romano procedure under the weakly dependent standard Gaussian setting with correlation matrix $\Sigma_n$, i.e.,
\[k\text{-}FWER_{BS}(n, \alpha, \Sigma_n) \longrightarrow \displaystyle1- e^{-k\alpha} \sum_{s=0}^{k-1} \frac{(k\alpha)^s}{s!} \text { as } n \rightarrow \infty.\]
\end{remark}
Following the one-sided case, we derive analogous asymptotic power results for the present two-sided testing problem. The following two theorems depict the asymptotics.

\begin{theorem}\label{bs_pwr_1}
    Suppose $n_1\to\infty$ and $\frac{n_1}{n}\to p_1\in(0,\,1]$ as $n\to\infty$. Then, under the weakly correlated Gaussian setting, for both \autoref{test_bon_bs} and \autoref{test_sidak_bs}, $\displaystyle\lim_{n\to\infty}\,AnyPwr_{BS}=1$ if $\displaystyle\lim_{n_1\to\infty}\,\frac{\sqrt{2\,\log n_1}}{\mu_{(n_1)}}<1$, where $\displaystyle\mu_{(n_1)}=\max_{i\in\mathcal{I_1}}|\mu_i|$.
\end{theorem}

\begin{theorem}\label{bs_pwr_2}
    Suppose $\displaystyle\lim_{n\to\infty}\min_{i\in\mathcal{I}_1}|\mu_i|=\mu>0$. Then, if $\displaystyle\lim_{n\to\infty}\frac{n_1}{n}e^{\mu\sqrt{2\log n}}=\infty$, $\displaystyle\lim_{n\to\infty}\,AnyPwr_{BS}=1$, for both \autoref{test_bon_bs} and \autoref{test_sidak_bs}.
\end{theorem}

\section{Simulation Studies \label{sec6}}

The FWER of a single-step multiple testing procedure (using the cut-off $c_n$ for each of the hypotheses) for the equicorrelated Normal setup under the global null is given by \citep{dey2025stpa}:
\begin{equation}
    FWER_{H_0}=1-\mathbb{E} \left[\Phi^n\left(\frac{c_n + \sqrt{\rho} Z}{\sqrt{1 - \rho}}\right)\right], \text{ where } Z \sim N(0,\,1). \label{eqFWERrho}
\end{equation}
Here $\rho$ is the common correlation. This essentially elucidates the ease of simulating from a equicorrelated Normal setup. Simulating Normal random variables having more general dependencies might be extremely computationally intensive (depending on the correlation structure) because of the presence of a large, and much more complex covariance matrix. Here we shall work with product correlation structures: 
$$\forall\ i\neq j, \quad \rho_{ij} = \lambda_i \lambda_j, \quad \text{where $\lambda_i \in (-1,\,1)$ for each $i$.}$$

\noindent \cite{tong2} refers to this as \textit{structure} $\mathbf{\ell}$. Now, suppose $X_1,\,X_2,\,\dots,\,X_n$ are standard Normal variables having $Corr(X_i,\,X_j)=\rho_{ij}=\lambda_i\lambda_j$ for all $i\ne j$, $\lambda_i\in(-1,\,1)$. Equivalently, one can write $X_i=\lambda_iZ+ \sqrt{1-\lambda_i^2}\ \epsilon_i$ for all $i=1,\,2,\,\dots,\,n$, where $Z, \epsilon_1,\,\epsilon_2,\,\dots,\,\epsilon_n$ are independent $N(0,\,1)$ variables.
One can generalize \eqref{eqFWERrho} to this case as follows.

\noindent Consider a single-step MTP using common right-sided cutoff $\tau_n$. The FWER of this procedure is:
\begin{align}\label{product_corr_1}
FWER_{H_0}
&= 1 - \mathbb{P}_{H_0} \left( \bigcap_{i=1}^n \{ X_i \leq \tau_n \} \right) \notag\\
&= 1 - \mathbb{P}_{H_0} \left( \bigcap_{i=1}^n \left\{ \lambda_iZ\ +\ \sqrt{1-\lambda_i^2}\ \epsilon_i \leq \tau_n \right\} \right) \notag\\
&=1 - \mathbb{P}_{H_0} \left( \bigcap_{i=1}^n \{\epsilon_i\leq \frac{\tau_n\ -\ \lambda_iZ}{\sqrt{1-\lambda_i^2}}\} \right)\notag\\
&=1 - \mathbb{E}_Z \left[ \prod_{i=1}^n \Phi \left( \frac{\tau_n - \lambda_i Z}{\sqrt{1 - \lambda_i^2}} \right) \right], \quad \text{where } Z \sim \mathcal{N}(0,\,1).
\end{align}
Similarly, for a single-step MTP using common both-sided cutoffs $-\tau_n$ and $\tau_n$, the FWER of this procedure is
\begin{equation}
    FWER_{H_0}=1 - \mathbb{E}_Z \left[ \prod_{i=1}^n \left\{\Phi \left( \frac{\tau_n - \lambda_i Z}{\sqrt{1 - \lambda_i^2}} \right)-\Phi \left( \frac{-\tau_n - \lambda_i Z}{\sqrt{1 - \lambda_i^2}} \right)\right\} \right], \quad \text{where } Z \sim \mathcal{N}(0,\,1). \label{product_corr_2}
\end{equation}

\subsection{Simulation Scheme}

For our purpose, we need to simulate FWER under global null for multivariate standard Normal distribution under weak dependence structure, \eqref{wd_cond}. For this, we utilize the notion of product correlation. So, we find real numbers $\lambda_1,\,\lambda_2,\,\dots,\,\lambda_n$, each in $(0,\,1)$ such that $\rho_{ij}\geq 0$ and the weak dependence condition \eqref{wd_cond} is satisfied. In particular, if we take 
\begin{equation}
\lambda_i =
\begin{cases}
\lambda & \text{for } i = 1, \\
\frac{1}{(\log i)^{1+\delta}} & \text{for } i = 2,\,3,\,\dots,\,n, \text{ for some } \lambda>0 \text{ and } \delta>0,
\end{cases}
\label{eq:lambda_def}
\end{equation}
then for this product correlation structure the weak dependency condition is satisfied. Hence, considering the product correlation structure defined in \eqref{eq:lambda_def}, we simulate FWER for both \autoref{test_bon} and \autoref{test_sidak} using \eqref{product_corr_1} in the following way :

For given $(n,\alpha)$, we compute $c_{Bon}(n,\,\alpha)=\Phi^{-1}(1-\frac{-\log(1-\alpha)}{n})$ and $c_{Sid}(n,\,\alpha)=\Phi^{-1}\left((1-\alpha)^{1/n}\right)$. Then, by choosing some $\lambda>0$ and $\delta>0$, we compute $\lambda_1,\,\lambda_2,\,\dots,\,\lambda_n$ according to \eqref{eq:lambda_def}. Now, we generate 10000 independent observations from $N(0,\,1)$ (these are the $Z$ variables) and for each of these observations, we correspondingly compute $\Phi \left( \frac{c_n(\alpha) - \lambda_i Z}{\sqrt{1 - \lambda_i^2}} \right)$ for all $n$ $\lambda_i$' s and take the product of them to obtain $\prod_{i=1}^{n} \Phi\left( \frac{c_n(\alpha)\ -\ \lambda_i Z}{\sqrt{1 - \lambda_i^2}} \right)$. Then taking the mean of this quantity over all 10000 $Z$ observations and subtracting this mean from 1, we obtain the simulated FWER under global null for this setting. Analogous steps are carried out for the both-sided case.

We repeat this whole process for several combinations of $(n, \alpha, \delta)$. Tables are provided which present the estimates of FWER (under $H_0$) for both the testing procedures, under a weakly dependent multivariate Normal setup satisfying the product correlation structure, as defined in previous section. For both the procedures, there are two tables for two choices of $\alpha$, viz., 0.1 and 0.05, considering different values of $n$ and $\delta$ in each table. 

\subsection{Simulation Results}
\subsubsection{Adjusted Bonferroni (One-sided)}

\renewcommand{\arraystretch}{1.2} 
\setlength{\tabcolsep}{10pt} 


\vspace{5pt}

\begin{table}[H]
\centering
\captionsetup{skip=6.5pt}
\caption{Results for \(\alpha=0.10\)}
\begin{tabularx}{\textwidth}{l|*{4}{>{\centering\arraybackslash}X|}>{\centering\arraybackslash}X}
\hline
\textbf{\(\)} & \(\delta=0.1\) & \(\delta=0.25\) & \(\delta=0.5\) & \(\delta=0.75\) & \(\delta=1\) \\
\hline
\(n=2500\)  & 0.09887 & 0.09883 & 0.09989 & 0.09985 & 0.09981 \\
\(n=5000\)  & 0.09765 & 0.09979 & 0.09981 & 0.09993 & 0.09980 \\
\(n=7500\)  & 0.09949 & 0.09938 & 0.09960 & 0.09991 & 0.09992 \\
\(n=10000\) & 0.09927 & 0.09892 & 0.10012 & 0.10007 & 0.10007 \\
\hline
\end{tabularx}
\end{table}

\vspace{5pt}
\begin{table}[H]
\centering
\captionsetup{skip=6.5pt}
\caption{Results for \(\alpha=0.05\)}
\begin{tabularx}{\textwidth}{l|*{4}{>{\centering\arraybackslash}X|}>{\centering\arraybackslash}X}
\hline
\textbf{\(\)} & \(\delta=0.1\) & \(\delta=0.25\) & \(\delta=0.5\) & \(\delta=0.75\) & \(\delta=1\) \\
\hline
\(n=2500\)  & 0.04971 & 0.04972 & 0.04979 & 0.05000 & 0.05003 \\
\(n=5000\)  & 0.04986 & 0.04989 & 0.05008 & 0.04984 & 0.05007 \\
\(n=7500\)  & 0.04928 & 0.05018 & 0.04994 & 0.04998 & 0.05003 \\
\(n=10000\) & 0.04983 & 0.04980 & 0.05011 & 0.05002 & 0.04999 \\
\hline
\end{tabularx}
\end{table}


\subsubsection{Sidak (One-sided)}
\renewcommand{\arraystretch}{1.2} 
\setlength{\tabcolsep}{10pt} 


\vspace{5pt}

\begin{table}[H]
\centering
\captionsetup{skip=6.5pt}
\caption{Results for \(\alpha=0.10\)}
\begin{tabularx}{\textwidth}{l|*{4}{>{\centering\arraybackslash}X|}>{\centering\arraybackslash}X}
\hline
\textbf{\(\)} & \(\delta=0.1\) & \(\delta=0.25\) & \(\delta=0.5\) & \(\delta=0.75\) & \(\delta=1\) \\
\hline
\(n=2500\)  & 0.09888 & 0.09884 & 0.09989 & 0.09986 & 0.09981 \\
\(n=5000\)  & 0.09765 & 0.09979 & 0.09981 & 0.09993 & 0.09980 \\
\(n=7500\)  & 0.09949 & 0.09938 & 0.09961 & 0.09992 & 0.09992 \\
\(n=10000\) & 0.09927 & 0.09892 & 0.10012 & 0.10007 & 0.10007 \\
\hline
\end{tabularx}
\end{table}

\vspace{5pt}
\begin{table}[H]
\centering
\captionsetup{skip=6.5pt}
\caption{Results for \(\alpha=0.05\)}
\begin{tabularx}{\textwidth}{l|*{4}{>{\centering\arraybackslash}X|}>{\centering\arraybackslash}X}
\hline
\textbf{\(\)} & \(\delta=0.1\) & \(\delta=0.25\) & \(\delta=0.5\) & \(\delta=0.75\) & \(\delta=1\) \\
\hline
\(n=2500\)  & 0.04971 & 0.04972 & 0.04979 & 0.05000 & 0.05003 \\
\(n=5000\)  & 0.04986 & 0.04989 & 0.05008 & 0.04984 & 0.05007 \\
\(n=7500\)  & 0.04928 & 0.05018 & 0.04994 & 0.04998 & 0.05003 \\
\(n=10000\) & 0.04983 & 0.04980 & 0.05011 & 0.05002 & 0.04999 \\
\hline
\end{tabularx}
\end{table}

\subsubsection{Adjusted Bonferroni (Two-sided)}

\renewcommand{\arraystretch}{1.2} 
\setlength{\tabcolsep}{10pt} 

\vspace{5pt}

\begin{table}[H]
\centering
\captionsetup{skip=6.5pt}
\caption{Results for \(\alpha=0.10\)}
\begin{tabularx}{\textwidth}{l|*{4}{>{\centering\arraybackslash}X|}>{\centering\arraybackslash}X}
\hline
\textbf{\(\)} & \(\delta=0.1\) & \(\delta=0.25\) & \(\delta=0.5\) & \(\delta=0.75\) & \(\delta=1\) \\
\hline
\(n=2500\)  & 0.09972 & 0.09989 & 0.09994 & 0.09995 & 0.09999 \\
\(n=5000\)  & 0.09980 & 0.10018 & 0.09996 & 0.10002 & 0.09998 \\
\(n=7500\)  & 0.10012 & 0.10001 & 0.09994 & 0.09997 & 0.09998 \\
\(n=10000\) & 0.10011 & 0.09991 & 0.10002 & 0.09999 & 0.10000 \\
\hline
\end{tabularx}
\end{table}

\vspace{5pt}
\begin{table}[H]
\centering
\captionsetup{skip=6.5pt}
\caption{Results for \(\alpha=0.05\)}
\begin{tabularx}{\textwidth}{l|*{4}{>{\centering\arraybackslash}X|}>{\centering\arraybackslash}X}
\hline
\textbf{\(\)} & \(\delta=0.1\) & \(\delta=0.25\) & \(\delta=0.5\) & \(\delta=0.75\) & \(\delta=1\) \\
\hline
\(n=2500\)  & 0.04993 & 0.04996 & 0.05002 & 0.05007 & 0.05003 \\
\(n=5000\)  & 0.04988 & 0.04999 & 0.04997 & 0.04999 & 0.05000 \\
\(n=7500\)  & 0.05004 & 0.05000 & 0.04997 & 0.05001 & 0.04999 \\
\(n=10000\) & 0.05004 & 0.04999 & 0.05000 & 0.04999 & 0.05000 \\
\hline
\end{tabularx}
\end{table}

\subsubsection{Sidak (Two-sided)}
\renewcommand{\arraystretch}{1.2} 
\setlength{\tabcolsep}{10pt} 

\vspace{5pt}

\begin{table}[H]
\centering
\captionsetup{skip=6.5pt}
\caption{Results for \(\alpha=0.10\)}
\begin{tabularx}{\textwidth}{l|*{4}{>{\centering\arraybackslash}X|}>{\centering\arraybackslash}X}
\hline
\textbf{\(\)} & \(\delta=0.1\) & \(\delta=0.25\) & \(\delta=0.5\) & \(\delta=0.75\) & \(\delta=1\) \\
\hline
\(n=2500\)  & 0.09972 & 0.09989 & 0.09994 & 0.09995 & 0.09999 \\
\(n=5000\)  & 0.09980 & 0.10018 & 0.09996 & 0.10002 & 0.09998 \\
\(n=7500\)  & 0.10012 & 0.10001 & 0.09994 & 0.09997 & 0.09998 \\
\(n=10000\) & 0.10011 & 0.09991 & 0.10002 & 0.09999 & 0.10000 \\
\hline
\end{tabularx}
\end{table}

\vspace{5pt}
\begin{table}[H]
\centering
\captionsetup{skip=6.5pt}
\caption{Results for \(\alpha=0.05\)}
\begin{tabularx}{\textwidth}{l|*{4}{>{\centering\arraybackslash}X|}>{\centering\arraybackslash}X}
\hline
\textbf{\(\)} & \(\delta=0.1\) & \(\delta=0.25\) & \(\delta=0.5\) & \(\delta=0.75\) & \(\delta=1\) \\
\hline
\(n=2500\)  & 0.04993 & 0.04996 & 0.05002 & 0.05007 & 0.05003 \\
\(n=5000\)  & 0.04988 & 0.04999 & 0.04997 & 0.05000 & 0.05000 \\
\(n=7500\)  & 0.05004 & 0.05000 & 0.04997 & 0.05001 & 0.04999 \\
\(n=10000\) & 0.05004 & 0.04999 & 0.05000 & 0.04999 & 0.05000 \\
\hline
\end{tabularx}
\end{table}


    The simulation results convincingly illustrate that under this type of weakly dependent correlation structure, for a large number of hypotheses, the FWER under global null is extremely close to the target values $\alpha$ for all $\delta>0$ for both procedures.

\section{Concluding Remarks} \label{sec7}

Explicit evaluation of FWER or $k$-FWER requires the knowledge of joint distribution of test statistics under null hypotheses. While this computation is simple under independence, this becomes intractable under dependence \citep{DeyBhandari2025gfwer}. We revisit here the classical testing problem of Normal means under weak dependence. We focus on an adjusted version of Bonferroni procedure and the Sidak procedure. We show that both of them control FWER at the desired level exactly as the number of hypotheses approaches infinity. Towards obtaining these results, we establish several new probabilistic asymptotic results that might be insightful in probability theory too. 

Several works have focused on asymptotic FWER in multiple testing \citep{DasBhandari2025, Proschan}. The premise of this work is much more general from the aspects of dependence structures, general configurations of hypotheses, and power considerations. 

There are several intriguing extensions worth exploring. One interesting problem would be to study the limits of stepwise procedures under such weakly dependent structures. Another challenge is to study these from a Bayesian point of view \citep{bogdan}.

\bibliography{references}

\begin{thebibliography}{24}
\providecommand{\natexlab}[1]{#1}
\providecommand{\url}[1]{\texttt{#1}}
\expandafter\ifx\csname urlstyle\endcsname\relax
  \providecommand{\doi}[1]{doi: #1}\else
  \providecommand{\doi}{doi: \begingroup \urlstyle{rm}\Url}\fi

\bibitem[Bogdan et~al.(2011)Bogdan, Chakrabarti, Frommlet, and Ghosh]{bogdan}
M.~Bogdan, A.~Chakrabarti, F.~Frommlet, and J.~K. Ghosh.
\newblock {Asymptotic Bayes-optimality under sparsity of some multiple testing procedures}.
\newblock \emph{The Annals of Statistics}, 39\penalty0 (3):\penalty0 1551 -- 1579, 2011.
\newblock URL \url{https://doi.org/10.1214/10-AOS869}.

\bibitem[Cramer(1946)]{Cramer}
H.~Cramer.
\newblock \emph{Mathematical Methods of Statistics}.
\newblock Princeton University Press, 1946.

\bibitem[Das and Bhandari(2025)]{DasBhandari2025}
N.~Das and S.~K. Bhandari.
\newblock {FWER} for normal distribution in nearly independent setup.
\newblock \emph{Statistics \& Probability Letters}, 219:\penalty0 110340, 2025.
\newblock URL \url{https://www.sciencedirect.com/science/article/pii/S0167715224003092}.

\bibitem[Datta and Dey(2026)]{DATTA2026110609}
S.~Datta and M.~Dey.
\newblock An asymptotically exact multiple testing procedure under dependence.
\newblock \emph{Statistics \& Probability Letters}, 230:\penalty0 110609, 2026.
\newblock URL \url{https://www.sciencedirect.com/science/article/pii/S0167715225002548}.

\bibitem[Dey(2024)]{deystpa}
M.~Dey.
\newblock On limiting behaviors of stepwise multiple testing procedures.
\newblock \emph{Statistical Papers}, 65:\penalty0 5691--5717, 2024.
\newblock URL \url{https://doi.org/10.1007/s00362-024-01613-6}.

\bibitem[Dey(2025)]{dey2025stpa}
M.~Dey.
\newblock Asymptotics in multiple hypotheses testing under dependence: beyond normality.
\newblock \emph{Statistical Papers}, 66\penalty0 (7), 2025.
\newblock URL \url{https://doi.org/10.1007/s00362-025-01770-2}.

\bibitem[Dey and Bhandari(2023)]{deybhandari}
M.~Dey and S.~K. Bhandari.
\newblock {FWER} goes to zero for correlated normal.
\newblock \emph{Statistics \& Probability Letters}, 193:\penalty0 109700, 2023.
\newblock URL \url{https://www.sciencedirect.com/science/article/pii/S0167715222002139}.

\bibitem[Dey and Bhandari(2024)]{deybhandaristpa}
M.~Dey and S.~K. Bhandari.
\newblock Bounds on generalized family-wise error rates for normal distributions.
\newblock \emph{Statistical Papers}, 65:\penalty0 2313–2326, 2024.
\newblock URL \url{https://doi.org/10.1007/s00362-023-01487-0}.

\bibitem[Dey and Bhandari(2025)]{DeyBhandari2025gfwer}
M.~Dey and S.~K. Bhandari.
\newblock Some results on generalized familywise error rate controlling procedures under dependence.
\newblock \emph{Stat}, 14\penalty0 (3):\penalty0 e70088, 2025.
\newblock URL \url{https://onlinelibrary.wiley.com/doi/abs/10.1002/sta4.70088}.

\bibitem[Dickhaus(2008)]{DickhausThesis}
T.~Dickhaus.
\newblock \emph{False Discovery Rate and Asymptotics, PhD Thesis}.
\newblock Heinrich-Heine-University Düsseldorf, 2008.

\bibitem[Dickhaus(2014)]{Dickhaus}
T.~Dickhaus.
\newblock \emph{Simultaneous Statistical Inference With Applications in the Life Sciences}.
\newblock Springer, Heidelberg, 2014.

\bibitem[Dickhaus et~al.(2026)Dickhaus, Heller, Hoang, and Rinott]{Dickhaus2026}
T.~Dickhaus, R.~Heller, A.~T. Hoang, and Y.~Rinott.
\newblock A procedure for multiple testing of partial conjunction hypotheses based on a hazard rate inequality.
\newblock \emph{Bernoulli}, 32\penalty0 (1):\penalty0 274--298, 2026.

\bibitem[Dudoit and Laan(2008)]{Dudoit}
S.~Dudoit and Mark~J. Laan.
\newblock \emph{Multiple Testing Procedures with Applications to Genomics}.
\newblock Springer Series in Statistics. Springer, 2008.

\bibitem[Farcomeni(2007)]{Farcomeni}
A.~Farcomeni.
\newblock Some results on the control of the false discovery rate under dependence.
\newblock \emph{Scandinavian Journal of Statistics}, 34\penalty0 (2):\penalty0 275--297, 2007.
\newblock URL \url{https://onlinelibrary.wiley.com/doi/abs/10.1111/j.1467-9469.2006.00530.x}.

\bibitem[Finner et~al.(2007)Finner, Dickhaus, and Roters]{FDR2007}
H.~Finner, T.~Dickhaus, and M.~Roters.
\newblock {Dependency and false discovery rate: Asymptotics}.
\newblock \emph{The Annals of Statistics}, 35\penalty0 (4):\penalty0 1432--1455, 2007.
\newblock URL \url{https://www.jstor.org/stable/25464546}.

\bibitem[Finner et~al.(2009)Finner, Dickhaus, and Roters]{FDR2009}
H.~Finner, T.~Dickhaus, and M.~Roters.
\newblock {On the false discovery rate and an asymptotically optimal rejection curve}.
\newblock \emph{The Annals of Statistics}, 37\penalty0 (2):\penalty0 596--618, 2009.
\newblock URL \url{https://doi.org/10.1214/07-AOS569}.

\bibitem[Gontscharuk and Finner(2013)]{GONTSCHARUK}
V.~Gontscharuk and H.~Finner.
\newblock Asymptotic fdr control under weak dependence: A counterexample.
\newblock \emph{Statistics \& Probability Letters}, 83\penalty0 (8):\penalty0 1888--1893, 2013.
\newblock URL \url{https://www.sciencedirect.com/science/article/pii/S0167715213001466}.

\bibitem[H{\"u}sler(1983)]{Husler1983}
J.~H{\"u}sler.
\newblock Asymptotic approximation of crossing probabilities of random sequences.
\newblock \emph{Zeitschrift f{\"u}r Wahrscheinlichkeitstheorie und Verwandte Gebiete}, 63\penalty0 (2):\penalty0 257--270, 1983.
\newblock URL \url{https://doi.org/10.1007/BF00538965}.

\bibitem[Leadbetter et~al.(1983)Leadbetter, Lindgren, and Rootzén]{Leadbetter}
M.~R. Leadbetter, G.~Lindgren, and H.~Rootzén.
\newblock \emph{Extremes and Related Properties of Random Sequences and Processes}.
\newblock Springer Series in Statistics, 1983.
\newblock URL \url{https://doi.org/10.1007/978-1-4612-5449-2}.

\bibitem[Lehmann and Romano(2005)]{Lehmann2005}
E.~L. Lehmann and J.~P. Romano.
\newblock {Generalizations of the familywise error rate}.
\newblock \emph{The Annals of Statistics}, 33\penalty0 (3):\penalty0 1138 -- 1154, 2005.
\newblock URL \url{https://doi.org/10.1214/009053605000000084}.

\bibitem[Proschan and Shaw(2011)]{Proschan}
M.~A. Proschan and P.~A. Shaw.
\newblock Asymptotics of {B}onferroni for dependent normal test statistics.
\newblock \emph{Statistics \& Probability Letters}, 81\penalty0 (7):\penalty0 739--748, 2011.
\newblock URL \url{https://www.sciencedirect.com/science/article/pii/S0167715210003329}.
\newblock Statistics in Biological and Medical Sciences.

\bibitem[Storey and Tibshirani(2003)]{StoreyTib}
J.~D. Storey and R.~Tibshirani.
\newblock Statistical significance for genomewide studies.
\newblock \emph{Proceedings of the National Academy of Sciences}, 100\penalty0 (16):\penalty0 9440--9445, 2003.
\newblock URL \url{https://www.pnas.org/doi/abs/10.1073/pnas.1530509100}.

\bibitem[Storey et~al.(2003)Storey, Taylor, and Siegmund]{Storey2004}
J.~D. Storey, J.~E. Taylor, and D.~Siegmund.
\newblock Strong control, conservative point estimation and simultaneous conservative consistency of false discovery rates: A unified approach.
\newblock \emph{Journal of the Royal Statistical Society Series B: Statistical Methodology}, 66\penalty0 (1):\penalty0 187--205, 2003.
\newblock URL \url{https://doi.org/10.1111/j.1467-9868.2004.00439.x}.

\bibitem[Tong(1980)]{tong2}
Y.~L. Tong.
\newblock \emph{Probability Inequalities in Multivariate Distributions}.
\newblock Probability and Mathematical Statistics. Academic Press, 1980.

\end{thebibliography}

\appendix

\section*{Appendix}

\section{Proofs of theoretical results mentioned in Section 3}

\subsection{Proof of \autoref{probabilistic_k_dn}}
We at first prove the result for $\tau\in(0,\,\infty)$ and then extend to the extreme cases.

Note that the correlation condition (1.1) in \cite{Husler1983} is equivalent to the correlation structure of the weakly dependent standard Gaussian sequence $\{X_n\}$.

Let $\{v_n\}$ be a sequence of $n$ such that $n\left(1-\Phi(v_n)\right)\to\tau\in(0,\,\infty)$ as $n\to\infty$. Thus, condition A of \cite{Husler1983} for the random sequence $\{X_n\}$ and the cutoff sequence $\{v_n\}$ clearly holds.

Let $N_n:=\# \{i\le n:X_i>v_n\}$. Then, utilizing \autoref{pwr_thm_1} and Theorem 3.5 of \cite{Husler1983}, we get
\[N_n \xrightarrow{d} \mathrm{Poisson}(\tau) \hspace{2mm} \text{as $n \to \infty$}.
\]
Therefore, for each fixed $k\ge 1$,
\begin{align}\label{eq_k_lar1}\displaystyle 
     \mathbb{P}\left(M_n^k\le v_n\right)=\mathbb{P}\left(N_n\le k-1\right)\notag
     &=\sum_{s=0}^{k-1}\mathbb{P}\left(N_n=s\right)\notag\\
     &\longrightarrow e^{-\tau} \sum_{s=0}^{k-1}\frac{\tau^s}{s!}\text{ as }n\to\infty.
 \end{align}
\noindent So, we are only required to show $\mathbb{P}\left(M_{d_n}^k\le u_n\right)-\mathbb{P}\left(M_{d_n}^k\le v_{d_n}\right)\to0\text{ as }n\to\infty$.
\begin{align}\label{eq_k_lar2}\displaystyle
    |\mathbb{P}\left(M_{d_n}^k\le u_n\right)-\mathbb{P}\left(M_{d_n}^k\le v_{d_n}\right)|&=\mathbb{P}\left(\min\{u_n, v_{d_n}\}\le M_{d_n}^k\le \max\{u_n, v_{d_n}\}\right)\notag\\
    &\le \mathbb{P}\left(\bigcup_{i=1}^{d_n}\{\min\{u_n, v_{d_n}\}\le X_i\le \max\{u_n, v_{d_n}\}\}\right)\notag\\
    &\le |d_n(\Phi(u_n)-\Phi(v_{d_n}))|\notag\\
    &=\bigg|d_n\left(1-\Phi(v_{d_n})\right)-d_n\left(1-\Phi(u_n)\right)\bigg|\notag\\
    &\longrightarrow0\text{ as }n\to\infty.\\
    &[\because \lim_{n\to\infty}d_n\left(1-\Phi(v_{d_n})\right)=\lim_{n\to\infty}d_n\left(1-\Phi(u_n)\right)=\theta\tau.]\notag
\end{align}
The extreme cases are remaining to be shown.

So, now we take $\{d_n\}$ and $\{u_n\}$ such that $d_n\left(1-\Phi(u_n)\right)\to0$ as $n\to\infty$.

Suppose, for each $\tau'\in(0,\,1)$, we have a sequence $\{w_n\}\left(w_n = w_n(\tau')\right)$ so that $d_n\left(1-\Phi(w_n)\right)\to\tau'$ as $n\to\infty$. Then, using \eqref{eq_k_lar1} and \eqref{eq_k_lar2} we have, $P\{M_{d_n}^k \le w_n\} \to e^{-\tau'}\sum_{s=0}^{k-1}\frac{{\tau'}^s}{s!}.$

Now if $\lim_{n \to \infty}d_n(1 - \Phi(u_n)) =0$, we must have $u_n > w_n$ for sufficiently large $n$, so that
\[
\limsup_{n\to\infty} P\{M_{d_n}^k \le u_n\}
\ge \lim_{n\to\infty} P\{M_{d_n}^k \le w_n\}
= e^{-\tau'}\sum_{s=0}^{k-1}\frac{{\tau'}^s}{s!}.
\]
Since this holds for arbitrary $\tau' > 0$, it follows that $P\{M_{d_n}^k \le u_n\} \to 1$ as $n\to\infty$.

The corresponding result for the case $\lim_{n \to \infty}d_n\left(1-\Phi(u_n)\right)=\infty$ follows similarly, completing the proof.

\subsection{Proof of \autoref{rate_of_convergence}}

We now analyze the convergence rates of the limiting results of \autoref{fwer_bon_2}. We consider the same re-notations and setup as mentioned in section 3 of the main paper. Let $c_{n}$ denote any of $c_{Bon}(n,\,\alpha)$ or $c_{Sid}(n,\,\alpha)$. Now, 
\begin{align}\label{rate_diff}
 |F W E R-\alpha|&=|\mathbb{P}\left(M_{n_{0}} \leq c_{n}\right)-(1-\alpha)| \notag\\
& \leq|\mathbb{P}\left(M_{n_{0}} \leq c_{n}\right)-\Phi^{n_{0}}\left(c_{n}\right)|+|\Phi^{n_{0}}\left(c_{n}\right)-(1-\alpha)| \notag \\
& =I_{1}+I_{2}  \quad (\text{say}).
\end{align}
Using Corollary 4.2.4 of \cite{Leadbetter}, 
$$
I_{1} \leq K \cdot \sum_{1 \leq i<j \leq n_{0}}|r_{i j}| \cdot \exp \left(-\frac{c_{n}^{2}}{1+|r_{i j}|}\right).
$$
Here $K$ depends on $\sup _{i, j} |r_{i j}|<1$, $r_{ij}=\operatorname{Corr}\left(Z_{i}, Z_{j}\right)$. Clearly
\begin{equation}\label{ineq_I1}
    I_{1} \leq K n \sum_{i =1}^n\rho_{i} \cdot \exp \left(-\frac{c_{n}^{2}}{1+\rho_i}\right).
\end{equation}
We know $\gamma=\sup_{n \geq 1}\rho_{n}<1$. Let $\nu$ be some constant such that $0<\nu<\frac{1-\gamma}{1+\gamma}$. Let 
$p=\left[n^{\nu}\right]$ and $\gamma_n = \sup_{m \geq n} \rho_{m}$.

\noindent The proof of Lemma 4.3.2 of \cite{Leadbetter} gives 
$$n \sum_{i=1}^n \rho_{i} \cdot \exp \left(-\frac{c^2_n}{1+\rho_{i}}\right) \leq K_{1} \cdot n ^{\frac{1+\nu-2}{1+\gamma}}\left(\log n\right)^{\frac{1}{1+\gamma}}+K_{2} \cdot \gamma_{p} c_n^{2} \exp \left(\gamma_{p} c_n^{2}\right).$$
Now, 
$$\gamma_{p} c_n^{2}=\gamma_{[n^{\nu}]} c_n^{2} \sim 2 \gamma_{[n^{\nu}]} \cdot \log n 
\sim \frac{2}{\nu} \gamma_{\left[n^{\nu}\right]} \log n^{\nu}.$$
Therefore, 
$$n \sum_{i=1}^{n} \rho_{i} \cdot \exp \left(-\frac{c^2_{n}}{1+\rho_{i}}\right) \leq K_{1}^{\star} \cdot n ^{\frac{1+\nu-2}{1+\gamma}}\left(\log n\right)^{\frac{1}{1+\gamma}}+K_{2}^{\star} \gamma_{\left[n^{\nu}\right]} \log n^{\nu}.$$
Using \eqref{ineq_I1} and the above inequality, we obtain
\begin{equation}\label{ineq_I1_final}
    I_{1} \leq l_{1} \cdot n^{\frac{1+\nu-2}{1+\gamma}}(\log n)^{\frac{1}{1+\gamma}}+l_{2} \cdot \gamma_{\left[n^{\nu}\right]} \cdot \log n^{\nu}.
\end{equation}
Now, 
\begin{align*}
    \Phi^{n_{0}}\left(c_{n}\right) & =\left(1-\frac{-\log (1-\alpha)-O\left(\frac{1}{n}\right)}{n}\right)^{n_{0}}
    =\exp \left[n_0 \cdot \log \left(1-\frac{-\log (1-\alpha)-O\left(\frac{1}{n}\right)}{n}\right)\right].
\end{align*}
One has $n_{0} \log \left(1-\frac{x}{n}\right)=\frac{n_{0}}{n} \cdot\left[-x-\frac{x^{2}}{2 n}-\frac{x^{3}}{3 n^{2}}+O\left(\frac{1}{n}\right)\right]
=\frac{n_{0}}{n}\left[-x+O\left(\frac{1}{n}\right)\right]$. Therefore, $\exp \left[n_{0} \log \left(1-\frac{x}{n}\right)\right]=\exp \left(-\frac{n_{0}}{n} x\right) \cdot \exp \left(O\left(\frac{1}{n}\right)\right) =e^{-\frac{n_{0}}{n} x} \cdot\left(1+O\left(\frac{1}{n}\right)\right)$. Hence, 
\begin{align*}
    \Phi^{n_{0}}\left(c_{n}\right)=e^{-\frac{n_{0}}{n} \cdot(-\log (1-\alpha))} \cdot e^{O\left(\frac{1}{n}\right)}\left(1+O\left(\frac{1}{n}\right)\right) \\
=(1-\alpha)^{\frac{n_{0}}{n}} \cdot\left(1+O\left(\frac{1}{n}\right)\right).
\end{align*}
Therefore, 
\begin{align}\label{ineq_I2}
I_2 & =|\Phi^{n_0}\left(c_n\right)-(1-\alpha)| \notag\\
& =(1-\alpha) \cdot \left|(1-\alpha)^{\frac{n_0}{n}-1} \left(1+O\left(\frac{1}{n}\right)\right)-1 \right| \notag\\
& = (1-\alpha)\cdot \left|\left\{1+\left(\frac{n_0}{n}-1\right)\log(1-\alpha)
+O\left(\left(\frac{n_0}{n}-1\right)^2\right)\right\}
\left(1+O\left(\frac{1}{n}\right)\right)-1\right|\notag\\
&=O\left(1-\frac{n_0}{n}\right)+O\left(\frac{1}{n}\right).
\end{align}

\noindent Inequalities \eqref{rate_diff}, \eqref{ineq_I1_final} and \eqref{ineq_I2} together imply
$$|FWER -\alpha| \leq l_{1} \cdot n^{\frac{1+\nu-2}{1+\gamma}}(\log n)^{\frac{1}{1+\gamma}}+l_{2} \cdot \gamma_{\left[n^\nu\right]} \cdot \log n^{\nu}
+O\left(1-\frac{n_0}{n}\right)+O\left(\frac{1}{n}\right).
$$

\section{Proofs of theoretical results mentioned in Section 4}

\subsection{Proof of \autoref{pwr_thm_1}}
   For any right-sided cutoff $\tau_n$, we have
\begin{align}\label{pwr_eq_1}
\ AnyPwr&=1-\mathbb{P}\left(\displaystyle\bigcap_{i=1}^{n_1}\{X_i\leq \tau_n\}\right)\notag\\
&\geq 1-\mathbb{P}(X_{n_1}\leq \tau_n)\quad[\text{where $X_{n_1}$ has mean $\mu_{n_1}$}]\notag\\
&=1-\mathbb{P}(Z\leq \tau_n-\mu_{n_1})\quad[Z\sim N(0,\,1)]\notag\\
&=1-\Phi(\tau_n-\mu_{n_1}).
\end{align}
\noindent \cite{Cramer} establishes, for $\beta>0$,
\begin{equation}\label{c_n_approx_1}\Phi^{-1}\left(1-\frac{\beta}{n}\right) = \sqrt{2\log n} - \frac{\log \log n + \log 4\pi+\log\beta}{2 \sqrt{2\log n}} + O \left(\frac{1}{\log n}\right).\end{equation}
\noindent Hence, using the expansion $(1-\alpha)^{\frac{1}{n}}=1-\frac{-\log(1-\alpha)-O(\frac{1}{n})}{n}$ and \eqref{c_n_approx_1}, we can write, for all sufficiently large $n$,
\begin{equation}\label{c_n_sid_1} \Phi^{-1}((1-\alpha)^{\frac{1}{n}})= \sqrt{2\log n} - \frac{\log \log n + \log 4\pi+\log(-\log(1-\alpha)-O(\frac{1}{n}))}{2 \sqrt{2\log n}} + O \left(\frac{1}{\log n}\right).\end{equation}
Now,
\begin{align*}\ c_{Bon}(n,\,\alpha)=\Phi^{-1}\left(1-\frac{-\ln(1-\alpha)}{n}\right)
&\sim \Phi^{-1}\left(1-\frac{-\ln(1-\alpha)p_1}{n_1}\right)\\
&\sim \sqrt{2\ln n_1}\quad(\text{using }\eqref{c_n_approx_1}).
\end{align*}
One also has
\begin{align*}\ c_{Sid}(n,\,\alpha)=\Phi^{-1}\left((1-\alpha)^{\frac{1}{n}}\right)
&\sim \Phi^{-1}\left([(1-\alpha)^{p_1}]^{\frac{1}{n_1}}\right)\\
&\sim \sqrt{2\ln n_1}\quad(\text{using }\eqref{c_n_sid_1}).
\end{align*}
Since $\lim_{n \to \infty}\frac{\sqrt{2\log n_1}}{\mu_{n_1}}<1$, the rest follows using \eqref{pwr_eq_1}.

\subsection{Proof of \autoref{Phi^dn}}
Using equation \eqref{c_n_approx_1},
\begin{equation}\label{c_n_approx_2}
    c_{\beta_n,\,n}^2  = 2\log n - \log \log n - \log 4\pi -\log\beta_n + o(1). \end{equation}
Hence, 
\begin{align}\label{c_n_prop}
    \frac{1}{\sqrt{2\pi}} \cdot \frac{ e^{-\frac{1}{2} c_{\beta_n,\,n}^2} \cdot e^{-c_{\beta_n,\,n}\cdot t}} {\sqrt{2 \log n}}
   \sim & \frac{e^{-\frac{1}{2}\left(2 \log n - \log \log n - \log 4 \pi-\log \beta_n \right)}}{\sqrt{2\pi}\cdot \sqrt{2 \log n}} \cdot e^{-t(c_{\beta_n,\,n} - \sqrt{2 \log n})}  \cdot e^{-t \sqrt{2 \log n}} \quad \text{(using \eqref{c_n_approx_2})}\notag\\
   \sim & \frac{\frac{1}{n} \cdot \sqrt{\log n}\cdot 2\sqrt{\pi}\cdot\sqrt{\beta}}{\sqrt{2\pi}\cdot \sqrt{2 \log n}} \cdot e^{-t \sqrt{2 \log n}} \quad \text{(since $c_{\beta_n,\,n}-\sqrt{2 \log n} \to 0$ from \eqref{c_n_approx_1})}\notag\\
   \sim & \frac{\sqrt{\beta}\cdot e^{-t \sqrt{2 \log n}}}{n}.
\end{align} 
Now, 
\begin{align}
d_n \cdot\left(1-\Phi\left(c_{\beta,\,n}+t\right)\right)\sim & d_n \cdot \frac{\phi\left(c_{\beta,\,n}+t\right)}{c_{\beta,\,n}+t} \notag\\
 \sim & d_n \cdot \frac{e^{-t^2/2}}{\sqrt{2\pi} } \cdot \frac{ e^{-\frac{1}{2} c_{\beta,\,n}^2} \cdot e^{- c_{\beta,\,n} t}}{\sqrt{2 \log n}} \notag\\
 = &  \sqrt{\beta} \cdot e^{-t^2/2} \cdot d_n \cdot \frac{e^{-t \sqrt{2 \log n}}}{n} \quad \text{(utilizing \eqref{c_n_prop}).}\notag
\end{align}
The rest is obvious.

\section{Proofs of theoretical results mentioned in Section 5}

The following results are imperatively necessary in proving \autoref{probabilistic_k_dn_bs}:
\begin{lemma}[Lemma 11.1.2 of \citet{Leadbetter}]\label{bs_lem1}
Let $\xi_1,,\,\xi_2,\,\dots,\,\xi_n$ be standard normal variables with covariance 
matrix $\Lambda^1 = (\Lambda^1_{ij})$ and $\eta_1,\,\eta_2,\,\dots,\,\eta_n$
be standard normal variables with covariance matrix $\Lambda^0 = (\Lambda^0_{ij})$. Let 
$\lambda_{ij} = \max(|\Lambda^1_{ij}|, |\Lambda^0_{ij}|)$.

Further, let $\mathbf{u} = (u_1,\,\dots,\,u_n)$ and 
$\mathbf{v} = (v_1,\,\dots,\,v_n)$ be vectors of real numbers and write \(w = \min\left(|u_1|,\,\dots,\,|u_n|,\,|v_1|,\,\dots,\,|v_n|\right).\)Then,
\begin{equation}
\begin{split}
\mathbb{P}\{-v_j < \xi_j \le u_j \text{ for } j=1,\,2,\,\dots,\,n\}
- \mathbb{P}\{-v_j < \eta_j \le u_j \text{ for } j=1,\,2,\,\dots,\,n\} \\
\le \frac{4}{2\pi}
\sum_{1 \le i < j \le n}
|\Lambda^1_{ij} - \Lambda^0_{ij}|
(1 - \lambda_{ij}^2)^{-1/2}
\exp\!\left(- \frac{w^2}{1 + \lambda_{ij}}\right).
\end{split}
\end{equation}

\end{lemma}

\begin{lemma}[Lemma 4.3.2 of \citet{Leadbetter}]\label{bs_lem2}

Suppose $r_n \log n \to 0$, and that $\{u_n\}$ is a sequence of constants such that $n(1 - \Phi(u_n))$ is bounded. Then,
\[n \sum_{j=1}^n |r_j|\,\exp\!\left(- \frac{u_n^2}{\,1 + |r_j|\,}\right)\;\longrightarrow\; 0\text{ as } n \to \infty .\]
\end{lemma}

\subsection{Proof of \autoref{probabilistic_k_dn_bs}}
We consider the case when $\tau\in(0,\,\infty)$ at first. Let $\Phi^{\star}$ be the common cdf of $|X_i|$'s. Hence, $\Phi^{\star}(u)=\Phi(u)-\Phi(-u)=2\Phi(u)-1$.

Let $\{v_n\}$ be a sequence of $n$ such that $n\left(1-\Phi(v_n)\right)\to\tau\in(0,\,\infty)$ as $n\to\infty$. So, $n\left(1-\Phi^{\star}(v_n)\right)=2n\left(1-\Phi(v_n)\right)\to2\tau\text{ as }n\to\infty$. This substantiates condition A of \cite{Husler1983} for the random sequence $\{|X_n|\}$ and the cutoff sequence $\{v_n\}$.
    
Now, let $1\le l_1< l_2< \cdots< l_p\le n$ for some fixed $p\le n$. Then, from \autoref{bs_lem1},
    \begin{align}\label{bs_eq1.1}
&\bigg|\mathbb{P}(|X_i|\le v_n, i=1,\,2,\,\dots,\,n)
- \left(\Phi^{\star}(v_n)\right)^p\bigg|\notag \\
&\le \frac{4}{2\pi}\sum_{i<j=l_1,\,\dots,\,l_p}\frac{|\rho_{ij}|}{\sqrt{(1-\rho_{ij}^2)}}\cdot\exp\!\left(- \frac{v_n}{1+|\rho_{ij}|}\right)\notag\\
&\le K\sum_{i<j=l_1,\,\dots,\,l_p}|\rho_{ij}|\cdot\exp\!\left(- \frac{v_n^2}{1+|\rho_{ij}|}\right) \quad [\text{where }K=\frac{4}{2\pi}\cdot\frac{1}{\sqrt{1-\gamma^2}}]\notag\\
&\le Kn\sum_{i=1}^n\rho_{i}\cdot\exp\!\left(- \frac{v_i^2}{1+\rho_{i}}\right).
\end{align}

Thus, for any $1\le i_1< i_2< \cdots< i_p<j_1<j_2<\cdots<j_q\le n$, for any fixed $p\ge 1, q\ge 1$ such that $p+q\le n$,
\begin{align}\label{cond_d}\bigg|\mathbb{P}(|X_i|\le v_n, i=i_1,\,\dots,\,i_p,\,j_1,\,\dots,\,j_q)&-\mathbb{P}(|X_i|\le v_n, i=i_1,\,\dots,\,i_p)\cdot\mathbb{P}(|X_i|\le u_n, i=j_1,\,.\dots,\,j_q)\bigg|\notag\\
&\le 3Kn\sum_{i=1}^n\rho_{i}\cdot\exp\!\left(- \frac{v_i^2}{1+\rho_{i}}\right)\notag\\
&\to0\text{ as }n\to\infty\text{ using \autoref{bs_lem2}}.
\end{align}

\noindent Hence, \eqref{cond_d} verifies condition D of \cite{Husler1983}.

\noindent For any $1\le i\ne j\le n$, utilizing \autoref{bs_lem1},
\begin{equation}\label{bs_eq1.2}\mathbb{P}\left(|X_i|> v_n, |X_j|> v_n\right)
\le \left(1-\Phi^{\star}(v_n)\right)^2+K\cdot|\rho_{ij}|\cdot\exp\!\left(- \frac{v_n}{1+|\rho_{ij}|}\right).\end{equation}
Now, let $r\le n$ be some integer and $I$ be a subset of $\{1,\,2,\,\dots,\,n\}$ such that $|I|:=\#\{i\in I\}\le n/r$. Then, using \eqref{bs_eq1.2} we have,
\begin{align}\label{cond_d'}
    \sum_{i<j\in I}\mathbb{P}\left(|X_i|> v_n, |X_j|> v_n\right)&\le \sum_{i<j\in I}\left(1-\Phi^{\star}(v_n)\right)^2+K\sum_{i<j\in I}|\rho_{ij}|\cdot\exp\left(- \frac{v_n}{1+|\rho_{ij}|}\right)\notag\\
    &\le 4(\frac{n}{r})^2\left(1-\Phi(v_n)\right)^2+K\cdot\frac{n}{r}\sum_{i=1}^n\rho_i\cdot\exp\left(- \frac{v_n}{1+\rho_i}\right)\notag\\
    &\to\frac{4\tau^2}{r^2}\text{ as }n\to\infty.
\end{align}

This implies $\lim_{r\to\infty}\lim_{n\to\infty}r\sum_{i<j\in I}\mathbb{P}\left(|X_i|> v_n, |X_j|> v_n\right)=0$, verifying condition D${^{\prime}}$ of \cite{Husler1983}.

Thus, by Theorem 3.5 of \cite{Husler1983}, letting $N_n^{\star}=\# \{i\le n:|X_i|>v_n\}$, we get $N_n \xrightarrow{d} \mathrm{Poisson}(2\tau)$ as $n \to \infty$.

Therefore, for each fixed $k\ge 1$,
\begin{align}\label{eq_k_lar_bs}\displaystyle 
     \mathbb{P}\left(L_n^k\le v_n\right)=\mathbb{P}\left(N_n^{\star}\le k-1\right)\notag
     &=\sum_{s=0}^{k-1}\mathbb{P}\left(N_n^{\star}=s\right)\notag\\
     &\longrightarrow e^{-2\tau} \sum_{s=0}^{k-1}\frac{(2\tau)^s}{s!}\text{ as }n\to\infty.
 \end{align}

\noindent So, we are only required to show $\mathbb{P}\left(L_{d_n}^k\le u_n\right)-\mathbb{P}\left(L_{d_n}^k\le v_{d_n}\right)\to0\text{ as }n\to\infty$.
\begin{align*}\displaystyle
    \left|\mathbb{P}\left(L_{d_n}^k\le u_n\right)-\mathbb{P}\left(L_{d_n}^k\le v_{d_n}\right)\right|&=\mathbb{P}\left(\min\{u_n, v_{d_n}\}\le L_{d_n}^k\le \max\{u_n, v_{d_n}\}\right)\\
    &\le \mathbb{P}\left(\bigcup_{i=1}^{d_n}\{\min\{u_n, v_{d_n}\}\le |X_i|\le \max\{u_n, v_{d_n}\}\}\right)\\
    &\le \big|d_n(\Phi^{\star}(u_n)-\Phi^{\star}(v_{d_n}))\big|\\
    &=2\bigg|d_n\left(1-\Phi(v_{d_n})\right)-d_n\left(1-\Phi(u_n)\right)\bigg|\\
    &\longrightarrow0\text{ as }n\to\infty.\\
    &[\because \lim_{n\to\infty}d_n\left(1-\Phi(v_{d_n})\right)=\lim_{n\to\infty}d_n\left(1-\Phi(u_n)\right)=\theta\tau.]
\end{align*}
Now, the result for the extreme cases, i.e., when $\lim_{n\to\infty}d_n\left(1-\Phi(u_n)\right)=0$ and $\lim_{n\to\infty}d_n\left(1-\Phi(u_n)\right)=\infty$ follows similarly as in the proof of \autoref{probabilistic_k_dn}. This completes the proof.

\subsection{Proof of \autoref{rate_of_convergence_bs}}
    $$
\begin{aligned}
& \mid FWER-\alpha \mid \\
= & \left| \mathbb{P}\left(|Z_{i}| \leq c_{2 n}, i=1(1) n_{0}\right) - (1-\alpha) \right|\\
\leq & \left|\mathbb{P}\left(|Z_{i}| \leq c_{2 n}, i=1(1) n_{0}\right) -\left(\Phi\left(c_{2 n}\right)-\Phi\left(-c_{2 n}\right)\right)^{n_{0}} \right|  +\left|\left(\Phi\left(c_{2 n}\right)-\Phi\left(-c_{2 n}\right)\right)^{n_{0}}-(1-\alpha)\right|.
\end{aligned}
$$

Now,
\begin{align}\label{bs_eq3}
    [\Phi(c_{2n}(\alpha))-\Phi(-c_{2n}(\alpha))]^{n_0}&=[2\Phi(c_{2d_n}(\alpha))-1]^{n_0}\notag\\
    &=[2\left(1-\frac{-\log(1-\alpha)-o(1)}{2n}\right)-1]^{n_0}\notag\\
    &=\left(1-\frac{-\log(1-\alpha)-o(1)}{n}\right)^{n_0}.
\end{align}
Using equations \eqref{bs_eq1.1} and \eqref{bs_eq3}, the rest follows along similar lines to the proof of \autoref{rate_of_convergence}.

\subsection{Proof of \autoref{bs_pwr_1}}
    For the both-sided testing problem, we have
    \begin{align*}
        AnyPwr_{BS}&=1-\mathbb{P}\left(\bigcap_{i\in\mathcal{I}_1}\{|X_i|\le c_{2n}(\alpha)\}\right)\\
        &=1-\mathbb{P}\left(\bigcap_{i\in\mathcal{I}_1}\{-c_{2n}(\alpha)-\mu_i\le Z_i\le c_{2n}(\alpha)-\mu_i\}\right), \quad Z_i=X_i-\mu_i\sim N(0,\,1),\\
        &\ge 1-\mathbb{P}\left(-c_{2n}(\alpha)-\mu_{(n_1)}\le Z\le c_{2n}(\alpha)-\mu_{(n_1)}\right), \quad Z\sim N(0,\,1),\\
        &=1-[\Phi\left(c_{2n}(\alpha)-\mu_{(n_1)}\right)-\Phi\left(-c_{2n}(\alpha)-\mu_{(n_1)}\right)]\\
        &\to 1\text{ as $n\to\infty$, under given condition.}
    \end{align*}

\subsection{Proof of \autoref{bs_pwr_2}}
    Let $\mathcal{I}_{11}$ and $\mathcal{I}_{12}$ denote the index sets with corresponding $X_i$'s having positive and negative means, respectively. Let
    $$U_i = \begin{cases} Z_i \quad &\text{for $i \in \mathcal{I}_{11}$,} \\
    -Z_i \quad &\text{for $i \in \mathcal{I}_{12}$.}\end{cases}$$
    Then, \begin{align*}
        & 1-AnyPwr_{BS}\\
        =& \mathbb{P}\left(\bigcap_{i\in\mathcal{I}_1}\{-c_{2n}(\alpha)-\mu_i\le Z_i\le c_{2n}(\alpha)-\mu_i\}\right), \quad Z_i=X_i-\mu_i\sim N(0,\,1),\\
        = &\mathbb{P}\left(\bigcap_{i\in\mathcal{I}_{11}}\{-c_{2n}(\alpha)-\mu_i\le Z_i\le c_{2n}(\alpha)-\mu_i\}\bigcap\bigcap_{j\in\mathcal{I}_{12}}\{-c_{2n}(\alpha)-\mu_j\le Z_j\le c_{2n}(\alpha)-\mu_j\}\right)\\
        \le &\mathbb{P}\left(\bigcap_{i\in\mathcal{I}_{11}}\{Z_i\le c_{2n}(\alpha)-\mu_i\}\bigcap\bigcap_{j\in\mathcal{I}_{12}}\{-c_{2n}(\alpha)-\mu_j\le Z_j\}\right)\\
        \le &\mathbb{P}\left(\bigcap_{i\in\mathcal{I}_{11}}\{Z_i\le c_{2n}(\alpha)-\mu\}\bigcap\bigcap_{j\in\mathcal{I}_{12}}\{-c_{2n}(\alpha)+\mu\le Z_j\}\right)\\
        =&\mathbb{P}\left(\bigcap_{i=1}^{n_1}\{U_i\le c_{2n}(\alpha)-\mu\}\right)\\
        \to & 0\text{ as }n\to\infty,\text{ using \autoref{Phi^dn} and \autoref{probabilistic_k_dn}, completing the proof.}
    \end{align*}

\end{document}